\documentclass[11pt]{article}

\usepackage{amsmath}
\usepackage{amssymb}
\usepackage{amscd}
\usepackage{amsthm}
\usepackage{indentfirst}
\usepackage[all]{xy}

\usepackage{hyperref}
\hypersetup{colorlinks=true,linkcolor=red,pdfpagemode=UseNone,pdfstartview=FitH}

\newcommand{\N}{\ensuremath{\mathbb{N}}}
\newcommand{\Z}{\ensuremath{\mathbb{Z}}}
\newcommand{\Q}{\ensuremath{\mathbb{Q}}}

\usepackage[margin=3cm]{geometry}

\theoremstyle{plain}
\newtheorem{theorem}{Theorem}[section]
\newtheorem{lemma}[theorem]{Lemma}
\newtheorem{corollary}[theorem]{Corollary}
\newtheorem{proposition}[theorem]{Proposition}
\theoremstyle{definition}
\newtheorem{definition}[theorem]{Definition}
\newtheorem{remark}[theorem]{Remark}
\newtheorem{example}[theorem]{Example}

\DeclareMathOperator{\Spec}{Spec}
\DeclareMathOperator{\Aut}{\underline{Aut}}
\DeclareMathOperator{\C}{\mathcal{C}}

\newcommand{\Site}{\mathcal{S}}
\newcommand{\fppf}{\text{\rm fppf}}
\newcommand{\kpl}{\text{\rm kfl}}
\newcommand{\gm}{\mathbf{G}_{{\rm m}}}

\DeclareMathOperator{\Hom}{\underline{Hom}}
\DeclareMathOperator{\GL}{GL}

\DeclareMathOperator{\Tors}{Tors}
\DeclareMathOperator{\Qcoh}{QCoh}
\DeclareMathOperator{\TC}{TameCover}

\begin{document}

\title{Tame stacks and log flat torsors}

\author{Jean Gillibert \and Heer Zhao}

\date{November 2023}

\maketitle

\begin{abstract}
We compare tame actions in the category of schemes with torsors in the category of log schemes endowed with the log flat topology. We prove that actions underlying log flat torsors are tame. Conversely, starting from a tame cover of a regular scheme that is an fppf torsor on the complement of a divisor with normal crossings, it is possible to build a unique log flat torsor that dominates this cover. In brief, the theory of log flat torsors gives a canonical approach to the problem of extending torsors into tame covers.
\end{abstract}



\section{Introduction}

Let $X$ be a noetherian scheme, and let $G$ be a finite flat group scheme over $X$. Let $Y$ be an $X$-scheme on which $G$ acts. What does it mean for the action of $G$ on $Y$ to be tame?

A large amount of work has been done on this question in the past few years. Let us first mention the foundational paper by Chinburg, Erez, Pappas and Taylor \cite{cept}, where a notion of tame action was introduced. This notion was used in many subsequent papers. In the present text, we refer to it as CEPT-tameness.

More recently, and independently, Abramovich, Olsson and Vistoli have introduced in \cite{aov} the notion of tame stack. Sophie Marques proved in \cite{sophie} that, if $X$ and $Y$ are affine and if the morphism $Y\rightarrow Y/G$ is flat, then these two notions of tameness coincide, namely, the action of $G$ on $Y$ is CEPT-tame if and only if the quotient stack $[Y/G]$ is tame (see Remark~\ref{sophie_remark}). Here $Y/G$ denotes the geometric and uniform categorical quotient in the category of algebraic spaces, which in the affine case is known to be a scheme. For more details, see Definition~\ref{def:tamestack} and the comments afterwards.

In this paper, the notion of tameness we use is that of \cite{aov}, but the applications we have in mind are linked with situations and problems relative to CEPT-tameness. Actually, our main results are formulated in a setting where the two notions coincide.

In fact, we consider in this paper a third point of view: that of Grothendieck and Murre \cite{gm71}, which is for us an important source of inspiration. The approach is to start from a $G$-torsor over some dense open subset of $X$, and try to extend it into some $X$-scheme on which $G$ acts tamely.

More precisely, let $X$ be a noetherian regular scheme, and let $D$ be a normal crossing divisor on $X$. Let $U\subseteq X$ be the open complement of $D$ in $X$. Roughly speaking, a tame $G$-cover of $X$ relative to $D$ is a finite flat morphism $Y\rightarrow X$ such that $G$ acts tamely on $Y$, and such that $Y_U\rightarrow U$ is a $G$-torsor (see Definition~\ref{tame_cover_def}).

The following questions immediately arise:

\begin{enumerate}
\item Given a $G$-torsor over $U$, under which condition is it possible to extend it into a tame $G$-cover of $X$?
\item Among all the possible extensions of a given $G$-torsor into a tame $G$-cover, is there a canonical choice?
\end{enumerate}

We give here an answer to these questions in terms of log flat torsors. More precisely, we show the following (cf. Theorem~\ref{main_thm}):

\begin{theorem}
\label{thm:main_intro}
Let $X$ be a noetherian regular scheme, and let $D$ be a normal crossing divisor on $X$, with complement $U:=X\setminus D$; we denote by $X(\log D)$ the log scheme defined by this data. Let $G$ be a finite flat group scheme over $X$, and let $t$ be a $G$-torsor over $U$. Then $t$ can be extended into a tame $G$-cover of $X$ if and only if it can be extended into a log flat $G$-torsor $t^{\log}$ over $X(\log D)$. If it exists, such a $t^{\log}$ is unique, and its underlying scheme together with the induced $G$-action is the initial object in the category of tame $G$-covers extending $t$.
\end{theorem}

When the group $G$ is \'etale, the first statement is implicitly contained in the work of Grothendieck and Murre \cite{gm71}: they require a tame cover to be a normal scheme, that is, they select the initial object right from the beginning. This implies the unicity of a tame cover extending a given torsor.

This work is a continuation of previous work of the first author \cite{gil6}, in which log schemes and log flat torsors were the central tool. The idea, already explored in \cite{gil6}, that tame actions should be linked with log flat torsors is now made precise here, at least in the particular setting of Theorem~\ref{thm:main_intro}. We also prove a more general result (see Corollary~\ref{vistoli_cor}) in a single direction: the action underlying a $G$-torsor for the log flat topology over any fs log scheme has a tame quotient stack. In \cite{gil6}, a similar result was proved under the hypothesis that $G$ is commutative. The main difference is that tame stacks are more flexible than CEPT-tame actions, that we were using previously.

Let us outline some consequences of Theorem~\ref{thm:main_intro}.

The theory of CEPT-tame actions and covers has been developed in order to extend to higher dimensional schemes methods and results about classical arithmetic Galois modules (\emph{e.g.} rings of integers of number fields), that are zero-dimensional objects over $\Spec(\Z)$. In section~\ref{CEPT_tameries}, we give examples where CEPT-tame covers induce log flat torsors. More precisely, given an arithmetic variety $X$ and a CEPT-tame $G$-cover $Y\rightarrow X$ that is flat and whose branch divisor has normal crossings, it is possible, according to Proposition~\ref{CEPT_prop}, to find another CEPT-tame $G$-cover $Y'\rightarrow X$ which is isomorphic to $Y$ outside the ramification locus, and can be endowed with the structure of a log flat $G$-torsor. We hope that this point of view will give a new insight into the work of Chinburg, Pappas, Taylor \emph{et al.}

Prior to that, another theory of tame actions, using the language of Hopf algebras, was introduced by Childs and Hurley \cite{ch86}, and subsequently developed by Childs. The initial aim was to turn wild extensions into tame objects under some Hopf algebra. This theory may nowdays be seen as a special case of CEPT-tameness. In section~\ref{CH_tameries}, we show that, over a Dedekind ring, if such a tame action is generically a torsor, then the maximal order on which the same Hopf algebra coacts can be endowed with the structure of a log flat torsor.

Let us review briefly the contents of this paper. In section~\ref{logflat_section} we prove the representability of log flat torsors under finite flat group schemes, which extends the result of Kato \cite{kato2} to the non-commutative case. Then we prove that, when $G$ is a finite flat linearly reductive group scheme, the restriction map from log flat torsors over $X(\log D)$ to fppf torsors over $U$ is an equivalence of categories ($X$, $D$ and $U$ being as in Theorem~\ref{thm:main_intro}). The main tools are a local description of $G$ as an extension of an \'etale group by a diagonalizable group, due to \cite{aov}, and a description of $G$-torsors similar to that given by Olsson in \cite{olsson}.

In section~\ref{tame_section}, we prove that a free action on a log scheme which is Kummer over the base has a tame quotient stack (in the category of schemes). In section~\ref{tame_divisor_section} we define tame $G$-covers relative to $D$ and we prove Theorem~\ref{main_thm}, which is the main result of this paper. The proof relies on the key point, proved in \cite{aov}, that tame stacks are \'etale-locally isomorphic to quotient stacks by finite flat linearly reductive group schemes. The end of the paper is devoted to applications, as outlined above.

\subsection*{Acknowledgements}

The first author expresses his warmest thanks to Angelo Vistoli for his help. Substantial improvements in this paper arose from discussions with him. He also thanks Dajano Tossici for inspiring conversations, and Luc Illusie and Ofer Gabber for stimulating exchanges. The first author acknowledges the support of the CIMI Labex. Part of the work had been done when the second author was a member of the research group of Professor Ulrich G\"ortz, partially under the support of Research Training Group 2553 (funded by the German Research Foundation DFG).


\section{Log flat torsors}
\label{logflat_section}

Unless otherwise stated, our log schemes are fine and saturated log schemes, and their log structures are defined on the small {\'e}tale site of their underlying scheme, which is assumed to be locally noetherian. In particular, finite flat schemes on our base scheme are the same as finite locally free schemes.

We endow the category of log schemes with the Kummer log flat topology, which was introduced by Kato \cite{kato2} under the name \emph{log flat topology}.


\subsection{Representability}
\label{rep_subsection}

In this section, we prove the representability of Kummer log flat torsors under finite flat group schemes. Our proof relies on the representability of $\mathrm{GL}_n$-torsors, proved by Kato \cite{kato2}.

\begin{theorem}
\label{thm_torsrep}
Let $X$ be an fs log scheme whose underlying scheme is locally noetherian. Let $G$ be a finite flat group scheme over the underlying scheme of $X$, endowed with the log structure induced by $X$.
Then any $G$-torsor in the category of sheaves for the Kummer log flat topology over $X$ is representable by an fs log scheme which is Kummer log flat over $X$, and whose underlying scheme is finite over the underlying scheme of $X$.
\end{theorem}

\begin{proof}
By a limit argument, we may assume that the underlying scheme of $X$ is $\Spec(A)$ for some noetherian strictly local ring $A$. Then $G$ is the spectrum of a free $A$-algebra of finite rank (say $n$), and the regular representation of $G$ yields a closed immersion $G\hookrightarrow \GL_{n,X}$.

Let $F$ be a $G$-torsor for the Kummer log flat topology on $X$. We have, as in the proof of \cite[Lemma~9.3]{kato2}, a cartesian square in the category of sheaves for the Kummer log flat topology on $X$
$$
\begin{CD}
\mathrm{GL}_n\times_X F @>>> \mathrm{GL}_n\wedge^G F \\
@V\mathrm{pr_1}VV @VVV \\
\mathrm{GL}_n @>>> G\backslash \mathrm{GL}_n
\end{CD}
$$
in which horizontal maps are obtained by moding out by the action of $G$, more precisely $\mathrm{GL}_n\wedge^G F$ is the quotient of $\mathrm{GL}_n\times_X F$ by the action defined by $g\cdot(a,b):=(ag,g^{-1}b)$, which is called the contracted product.

The horizontal map at the bottom is a classical fppf $G$-torsor, the two schemes being endowed with the log structures induced by $X$. The contracted product $\GL_n\wedge^G F$ is a $\GL_n$-torsor over $X$, hence is representable according to \cite[Prop.~9.7]{kato2}; it follows that the vertical map on the right is a morphism between log schemes. Therefore, the fiber product $\mathrm{GL}_n\times_X F$ is representable by an fs log scheme.

One deduces the representability of $F$ by considering the following cartesian square
$$
\begin{CD}
F @>1\times \mathrm{id}_F>> \mathrm{GL}_n\times_X F \\
@VVV @VV\mathrm{pr_1}V \\
X @>1>> \mathrm{GL}_n
\end{CD}
$$
where $1$ is the unit section of $\mathrm{GL}_n$.
Finally, the properties of $F\to X$ follow from Kummer log flat descent \cite[Theorem*~7.1]{kato2}.
\end{proof}

\begin{remark}
The representability of Kummer log flat $G$-torsors has been proved in \cite{kato2} and \cite{niziol} in the following cases:
\begin{enumerate}
\item $G$ is commutative, finite flat;
\item $G$ is commutative, smooth and affine;
\item $G=\GL_n$.
\end{enumerate}
Our Theorem~\ref{thm_torsrep} extends this result to any finite flat group scheme.
\end{remark}

We now prove that, under certain assumptions on the log structure, the morphism of schemes underlying a Kummer log flat torsor is classically flat.

\begin{proposition}
\label{prop:flatness}
Let $X$ be as in Theorem~\ref{thm_torsrep}, and let $G$ be a finite flat group scheme over the underlying scheme of $X$. We assume that for any geometric point $\bar{x}$ of $X$ the monoid $\overline{M}_{X,\bar{x}}$ is free. Let $f:Y\to X$ be a Kummer log flat $G$-torsor over $X$. Then the underlying map of schemes of $f$ is flat.
\end{proposition}

\begin{proof}
By Kummer log flat descent, see \cite[Theorem*~7.1]{kato2}, the map $f$ is Kummer log flat. By \cite[Chap. IV, Thm. 4.3.5]{ogu1}, it suffices to show that $f$ is integral (see \cite[Chap. III, Def. 2.5.1]{ogu1}). Let $y$ be any point of $Y$ and let $\bar{y}$ be any geometric point lying over $y$, we need to show that the homomorphism $\theta:\overline{M}_{X,f(\bar{y})}\to \overline{M}_{Y,\bar{y}}$ is an integral homomorphism of monoids. By \cite[Chap.~I, Prop.~4.7.5]{ogu1}, it suffices to show that $\theta$ is $\Q$-integral. Since $\theta$ is a Kummer homomorphism of finitely generated monoids, the $\Q$-integrality follows from the equivalence $1\Leftrightarrow 4$ from \cite[Chap. I, Thm. 4.7.7]{ogu1}.
\end{proof}

\begin{remark}
The statement above was proved in \cite[Prop.~3.2]{gil6} under the assumption that $G$ is commutative, which was there to ensure the representability of Kummer log flat $G$-torsors.
\end{remark}



\subsection{An equivalence of categories}

We now briefly recall a nice equivalence of categories that is due to Olsson \cite{olsson}.

Let $\Site$ be a site, with final object $e$. If $G$ is a group object in $\Site$, we denote by $\Tors_{\Site}(e,G)$ the category of $G$-torsors over $e$ in the category of sheaves on $\Site$.

Let 
$$
\begin{CD}
1 @>>> \Delta @>>> G @>s>> H @>>> 1 \\
\end{CD}
$$
be an exact sequence of group objects in $\Site$, where $\Delta$ is commutative. This implies that the map $G\rightarrow H$ is a $(\Delta,\Delta)$-bitorsor.

Let $P\rightarrow T$ be any $\Delta$-torsor. Then the contracted product $G\wedge^{\Delta} P$ is a $\Delta$-torsor over $H\times T$, the action being induced by the left action of $\Delta$ on $P$.

Let $P$ be a $G$-torsor. Let $\overline{P}$ be the quotient of $P$ by the action of $\Delta$. Then $P\rightarrow \overline{P}$ is a $\Delta$-torsor, $\overline{P}\rightarrow e$ is an $H$-torsor, and the action of $G$ on $P$ induces a map
\begin{align*}
\chi_P: G\wedge^{\Delta} P \;\longrightarrow & \; H\times P \\
(g,x) \;\longmapsto & \; (s(g),gx)
\end{align*}
which is a morphism between $\Delta$-torsors over $H\times \overline{P}$, and therefore is an isomorphism.
Putting all these data together allows us to define a category.

\begin{definition}
Let $\C$ be the category whose objects are triples
$$
(T,P\rightarrow T,\chi)
$$
where
\begin{enumerate}
\item $T\rightarrow e$ is an $H$-torsor
\item $P\rightarrow T$ is a $\Delta$-torsor
\item $\chi:G\wedge^{\Delta} P\rightarrow H\times P$ is a morphism of $\Delta$-torsors (over $H\times T$)
\item The diagram
$$
\begin{CD}
(G\wedge^{\Delta} G)\wedge^{\Delta} P @>m\wedge \mathrm{id}>> G\wedge^{\Delta} P \\
@| @| \\
G\wedge^{\Delta} (G\wedge^{\Delta} P) @>\mathrm{id}\wedge\chi>> G\wedge^{\Delta} P \\
\end{CD}
$$
commutes, where $m$ is the map induced by the multiplication in $G$.
\end{enumerate}
\end{definition}

\begin{theorem}
\label{olsson_torsors}
The functor
\begin{align*}
\Tors_{\Site}(e,G)\;\longrightarrow & \; \C \\
P\;\longmapsto & \; (\overline{P},P\rightarrow \overline{P},\chi_P)
\end{align*}
is an equivalence of categories.
\end{theorem}

\begin{proof}
The proof is extracted from \cite{olsson}. We will construct a quasi-inverse to this functor. Let $(T,P\rightarrow T,\chi)$ be an object of $\C$. Then by composing the maps
$$
\begin{CD}
G\times P @>>> G\wedge^{\Delta} P @>\chi >> H\times P @>\mathrm{pr}_2>> P \\
\end{CD}
$$
(where $\mathrm{pr}_2$ is the projection on the second factor) we get an action of $G$ on $P$. It is easy to check that this action endows $P$ with the structure of a $G$-torsor.
\end{proof}



\subsection{Extending torsors under linearly reductive group schemes}
\label{XlogD_subsection}

Let $X$ be a noetherian regular scheme, and let $D$ be a normal crossing divisor on $X$. Let $j:U\rightarrow X$ be the open complement of $D$ in $X$. We endow $X$ with the log structure defined by
$$
\mathcal{O}_X\cap j_*\mathcal{O}_U^*\longrightarrow \mathcal{O}_X.
$$
This log structure is just the direct image of the trivial log structure on $U$.
We denote by $X(\log D)$ the resulting log scheme, which is a fine and saturated log scheme.

In the case when $X$ is the spectrum of a discrete valuation ring and $D$ is the special point of $X$, the corresponding log structure is called the standard log structure on $X$.

If $G$ is an $X$-group scheme, we denote by $H^1_{\kpl}(X(\log D),G)$ the set of isomorphism classes of $G$-torsors for the Kummer log flat topology over $X(\log D)$, and by $H^1_{\fppf}(U,G_U)$ the set of isomorphism classes of $G_U$-torsors for the fppf topology over $U$.

The natural inclusion $U\rightarrow X(\log D)$ induces a localization morphism from the fppf site of $U$ to the Kummer log flat site of $X(\log D)$. This yields a natural restriction map
$$
\begin{CD}
H^1_{\kpl}(X(\log D),G) @>>> H^1_{\fppf}(U,G_U). \\
\end{CD}
$$
According to Proposition~\ref{prop:fullfaith}, this map is injective when $G$ is finite flat.
According to \cite[Prop.~3.2.1]{gil5}, this map is bijective when $G$ is a finite diagonalizable group scheme.

We shall prove that this map is bijective when $G$ is a finite flat linearly reductive group scheme over $X$. For that purpose, we first recall the definition given in \cite{aov}.

\begin{definition}
Let $S$ be a scheme, and let $G$ be an $S$-group scheme. We say that $G$ is linearly reductive if the functor
\begin{align*}
\Qcoh^G(S) &\longrightarrow \Qcoh(S) \\
F &\longmapsto F^G
\end{align*}
is exact.
\end{definition}

\begin{theorem}
\label{equivalence1}
Let $X$ be a noetherian regular scheme, let $D$ be a normal crossing divisor on $X$, and let $G$ be a finite flat linearly reductive group scheme over $X$. Then the restriction functor
$$
\begin{CD}
\Tors_{\kpl}(X(\log D),G) @>>> \Tors_{\fppf}(U,G_U) \\
\end{CD}
$$
is an equivalence of categories.
\end{theorem}

We first need a Lemma on Kummer log \'etale torsors.

\begin{lemma}
\label{log_etale_torsor}
Let $X$ and $D$ be as before, and let $f:Y\rightarrow X(\log D)$ be a Kummer log \'etale torsor under a finite constant group $H$. Then the underlying scheme of $Y$ is regular. Moreover, the complement of $V:=f^{-1}(U)$ is a normal crossing divisor on $Y$, and the log structure of $Y$ is associated to this divisor.
\end{lemma}

\begin{proof}
Let us first note (see \cite[Ex.~4.7 (c)]{illusie}) that $Y\rightarrow X(\log D)$ is a Kummer log \'etale $H$-torsor if and only if $V\rightarrow U$ is an \'etale $H$-torsor and the underlying scheme of $Y$ is a tamely ramified covering of $X$ relative to $D$ in the terminology of Grothendieck and Murre \cite[Def.~2.2.2]{gm71}. In this setting, the fact that the underlying scheme of $Y$ is regular was proved by Grothendieck and Murre \cite[Thm.~2.3.2 and Prop.~1.8.5]{gm71}.
Also, the complement of $V$ is a normal crossing divisor according to \cite[Lemma~1.8.6]{gm71}. The fact that the log structure of $Y$ is given by $V$ can be checked locally for the \'etale topology.
\end{proof}

\begin{proof}[Proof of Theorem~\ref{equivalence1}]
According to Proposition~\ref{prop:fullfaith} the restriction functor is fully faithful. When $G$ is a finite diagonalizable group scheme, it follows from \cite[Prop.~3.2.1]{gil5} that it is essentially surjective, hence it is an equivalence of categories. Let us prove the general case.

By descent of torsors, the question is local for the \'etale topology on $X$. Therefore, we may (and do) assume that $X$ is the spectrum of a strictly henselian local ring. In this context, it has been proved by \cite[Lemma~2.20]{aov} that we have an exact sequence
$$
\begin{CD}
1 @>>> \Delta @>>> G @>>> H @>>> 1 \\
\end{CD}
$$
where $\Delta$ is diagonalizable, and $H$ is constant of order coprime to the residue characteristics of $X$.
From this exact sequence, Theorem~\ref{olsson_torsors} gives us a description of $G$-torsors as triples. This description is valid in the topos of sheaves for the Kummer log flat topology on $X(\log D)$, as well as in the topos of sheaves for the fppf topology on $U$. So it remains to prove that the restriction functor
$$
(T,P\rightarrow T,\chi)\longmapsto (T_U,P_U\rightarrow T_U,\chi_U)
$$
is an equivalence between the two categories of triples. We will give a quasi-inverse to this functor. Let us consider a triple $(T_U,P_U\rightarrow T_U,\chi_U)$ over $U$. Since $H$ is constant of order coprime to the residue characteristics of $X$, there exists a unique Kummer log \'etale $H$-torsor $T\rightarrow X(\log D)$ extending $T_U$ \cite[Thm.~7.6]{illusie}. Moreover, by Lemma~\ref{log_etale_torsor}, we know that $T$ is a regular scheme, and that its log structure is induced by a normal crossing divisor. As pointed at the beginning of the proof, Theorem~\ref{equivalence1} holds for diagonalizable group schemes, hence (applying it over $T$) the $\Delta$-torsor $P_U\rightarrow T_U$ extends uniquely into a Kummer log flat $\Delta$-torsor $P\rightarrow T$. It remains to extend the map $\chi_U$. In other words, we have to prove that the restriction map
$$
\begin{CD}
\Tors_{\kpl}(H\times T,\Delta)@>>> \Tors_{\fppf}(H_U\times T_U,\Delta) \\
\end{CD}
$$
is fully faithful. But in fact this is an equivalence of categories, because $H$ being \'etale, the log scheme $T\times H$ is again a regular scheme whose log structure is associated to a divisor with normal crossings, and Theorem~\ref{equivalence1} holds for diagonalizable group schemes.
\end{proof}

\begin{remark}
Of course, Theorem~\ref{equivalence1} does not hold for an arbitrary linearly reductive group scheme. For example, the map
$$
\begin{CD}
H^1_{\kpl}(X(\log D),\gm) @>>> H^1_{\fppf}(U,\gm) \\
\end{CD}
$$
is not injective (see the description given in \cite{gil5}).
\end{remark}


\section{Tame stacks and free actions on log schemes}
\label{tame_section}


In this section, $X$ is a locally noetherian scheme, $(X,M)$ is a fine saturated log scheme, and $G$ is a finite flat group scheme over $X$. Given an $X$-scheme $Y$, an action of $G$ on $Y$ is a map $G\times_X Y\to Y$ satisfying the usual properties (associativity, unitarity). We say that the action is free if for all $T\to X$, the action of the group $G(T)$ on the set $Y(T)$ is free, \emph{i.e.} the induced map $G(T)\times Y(T) \to Y(T)\times Y(T);(g,y)\mapsto (gy,y)$ is injective. These definitions being functorial, they also make sense in the category of log schemes.

Let $Y\to X$ be an $X$-scheme on which $G$ acts. Let us recall the definition of the inertia group scheme of this action at some point of $Y$. Let $T$ be any scheme, and let $z:T\rightarrow Y$ be some element of $Y(T)$. Then the inertia group of $z$ is the fiber product
$$
I(z):=(Y\times G)\times_{Y\times Y} T
$$
where the product over $Y\times Y$ is relative to the maps $(m, \mathrm{pr}_1):Y\times G\rightarrow Y\times Y$ ($m$ being the action) and $(z,z):T\rightarrow Y\times Y$. It is easy to check that $I(z)$ is a subgroup scheme of $G_T$.

For the reader's convenience, we define the notion of tame stack in the special case of a quotient stack by a finite flat group scheme. We refer to \cite[Def.~3.1]{aov} for the general definition and properties of tame stacks. Everywhere in this paper, our stacks are algebraic stacks (in the sense of Artin).

\begin{definition}
\label{def:tamestack}
Let $f:Y\rightarrow X$ be an $X$-scheme locally of finite type on which $G$ acts, and let $q:Y\to Y/G$ be the geometric and uniform categorical quotient of this action in the category of algebraic spaces. The stack $[Y/G]$ is tame if the functor
$$
\begin{CD}
q_*^G:\Qcoh^G(Y) @>>> \Qcoh(Y/G) \\
\end{CD}
$$
is exact.
\end{definition}

The group scheme $G$ being finite flat over a locally noetherian base, the Keel-Mori theorem \cite{K-M} ensures us the existence of a geometric and uniform categorical quotient $Y/G$ in the category of algebraic spaces, which is the coarse moduli space of the stack $[Y/G]$. We refer to Conrad \cite{conrad2005} for a stacky version of the Keel-Mori theorem. Moreover, if the action of $G$ on $Y$ is admissible, i.e. if every $G$-orbit is contained in an open affine subscheme of $Y$, then $Y/G$ is actually a scheme.  In particular, if the $G$-invariant map $f:Y\to X$ is affine, then the action is admissible and $Y/G=\Spec_X((f_*\mathcal{O}_Y)^G)$. This condition will always be satisfied for the actions considered in this paper.

\begin{proposition}
\label{vistoli_prop}
Let $f:(Y,N)\rightarrow (X,M)$ be an $(X,M)$-log scheme on which $G$ acts. If the action of $G$ on $(Y,N)$ is free and $f$ is Kummer, then the inertia groups of the action of $G$ on $Y$ at geometric points of $Y$ are finite diagonalizable group schemes. In particular, the quotient stack $[Y/G]$ is tame.
\end{proposition}

\begin{proof}
According to \cite[Theorem~3.2]{aov}, in order to prove that $[Y/G]$ is tame, it suffices to show that, given a geometric point $x:\Spec(k)\rightarrow X$ of $X$, then for any $y:\Spec(k)\rightarrow Y$ above $x$, the inertia subgroup $I(y)\subseteq G_k$ is finite flat linearly reductive in the sense of \cite[Def.~2.2]{aov}. Finite diagonalizable group schemes being linearly reductive, the last part of the statement follows from the first one.

Pulling back $M$ along $x$ (resp. $N$ along $y$), we get two log structures $M_x$ and $N_y$ on $\Spec(k)$. The morphism $(Y,N)\rightarrow (X,M)$ induces a morphism $(\Spec(k),N_y)\rightarrow (\Spec(k),M_x)$. Now, let $\Aut((\Spec(k),N_y)\rightarrow (\Spec(k),M_x))$ be the group of automorphisms of the log scheme $(\Spec(k),N_y)$ that leave $(\Spec(k),M_x)$ invariant. Obviously, such automorphisms induce the identity on $\Spec(k)$. So what remains is the sheaf of automorphisms of the log structure $N_y$ that fix $M_x$. The field $k$ being algebraically closed, there exists a fine, saturated, sharp monoid $P$ (resp. $Q$) such that $M_x=k^*\oplus P$ (resp. $N_y=k^*\oplus Q$). Let $h:P=M_x/k^*\to N_y/k^*=Q$ be the canonical homomorphism. Then $h$ is a Kummer map of sharp fs monoids. We claim that $F:=\Aut((\Spec(k),N_y)\rightarrow (\Spec(k),M_x))=D({\rm coker}(h^{\mathrm{gp}}))$. Indeed let $\alpha:N_y\to N_y$ be an element of $F(\Spec(k))$. Then $\alpha$ is uniquely determined by three maps
\[k^*\xrightarrow{\alpha_{1}}k^*\oplus Q,Q\xrightarrow{\alpha_{21}}k^*,Q\xrightarrow{\alpha_{22}}Q.\]
Since $\alpha$ leaves $k^*\subset N_y$ invariant, $\alpha_1$ factors through $1_{k^*}$. We shall now prove that $\alpha_{22}=1_Q$. Since $\alpha$ leaves $f$ invariant, we have a commutative diagram
\[\xymatrix{
&P\ar@{^(->}[rd]^h\ar@{^(->}[ld]_h \\
Q\ar[rr]^{\alpha_{22}} &&Q
}.\]
Since $P$ and $Q$ are sharp and fs, they can be embedded in their groups of fractions $P^{\mathrm{gp}}$ and $Q^{\mathrm{gp}}$, which are finitely generated and torsion-free.
Since $h$ is Kummer, $h^{\mathrm{gp}}:P^{\mathrm{gp}}\to Q^{\mathrm{gp}}$ is injective and its cokernel is torsion, hence is finite (since $Q^{\mathrm{gp}}$ is finitely generated). So, the image of $h^{\mathrm{gp}}$ is a finite index subgroup of $Q^{\mathrm{gp}}$, and, by commutativity of the above diagram, the restriction of $\alpha_{22}^{\mathrm{gp}}$ to this subgroup is the identity. Since $Q^{\mathrm{gp}}$ is torsion-free, it follows that $\alpha_{22}^{\mathrm{gp}}=1_{Q^{\rm gp}}$, and thus $\alpha_{22}=1_Q$. So, the data of $\alpha$ is equivalent to the data of $\alpha_{21}$.

Since the composition $P\xrightarrow{h}Q\xrightarrow{\alpha_{21}}k^*$ is trivial, the data of $\alpha_{21}$ amounts to the data of a homomorphism ${\rm coker}(h^{\rm gp})\to k^*$. Since this argument is functorial under base-change by any morphism of schemes $T\to \Spec(k)$, we deduce that $F= \Hom({\rm coker}(h^{\rm gp}),\gm) = D({\rm coker}(h^{\mathrm{gp}}))$, and this finishes the proof of the claim.

Now, the action of $I(y)$ on $(\Spec(k),N_y)$ yields a map
$$
I(y)\longrightarrow D({\rm coker}(h^{\mathrm{gp}}))
$$
which is a monomorphism since the action of $G$ on $(Y,N)$ is free. Finally, we note that $I(y)$ is a closed subgroup scheme of the finite flat $k$-group scheme $G_k$, hence is finite flat. It follows that the map $I(y)\to D({\rm coker}(h^{\mathrm{gp}}))$ is a proper monomorphism, hence is a closed immersion by \cite[8.11.5]{EGA4.3}.
So $I(y)$ is a finite flat subgroup scheme of the diagonalizable group scheme $D({\rm coker}(h^{\mathrm{gp}}))$, hence is diagonalizable, as required.
\end{proof}

\begin{remark}
According to Proposition~\ref{vistoli_prop}, a free action on a log scheme which is Kummer over the base induces a tame action on the underlying schemes. These actions are quite special among tame actions: their inertia groups are finite diagonalizable group schemes, in particular they are commutative.
\end{remark}

We are now interested in $G$-torsors for the Kummer log flat topology. The statement $(i)$ below can be viewed as a generalization to the non-commutative case of \cite[Thm.~3.16]{gil6}, and statement $(ii)$ is some kind of Abhyankar's Lemma for Kummer log flat torsors. Let us stress here the fact that $(ii)$ is a consequence of the highly technical results of \cite{aov}.

\begin{corollary}
\label{vistoli_cor}
Let $f:(Y,N)\rightarrow (X,M)$ be a $G$-torsor for the Kummer log flat topology. Then:
\begin{enumerate}
\item[$(i)$] The quotient stack $[Y/G]$ is tame, with moduli space $X$.
\item[$(ii)$] Let $x\in X$ be a point. Then there exists an \'etale morphism $W\rightarrow X$ having $x$ in its image, a diagonalizable group scheme $\Delta$, and a finite scheme $T\rightarrow W$ on which $\Delta$ acts, together with an isomorphism $[Y/G]\times_X W\simeq [T/\Delta]$ of algebraic stacks over $W$.
\end{enumerate}
\end{corollary}

\begin{proof}
$(i)$ According to Proposition~\ref{vistoli_prop}, $[Y/G]$ is a tame stack. It remains to prove that its moduli space is $X$. Let us note that the action is admissible since $f:Y\rightarrow X$ is a finite (hence affine) $G$-invariant map by Theorem~\ref{thm_torsrep}. According to the Keel-Mori theorem \cite{K-M}, the moduli space of $[Y/G]$ is $\Spec_X((f_*\mathcal{O}_Y)^G)$, so it suffices to prove that $(f_*\mathcal{O}_Y)^G=\mathcal{O}_X$. Since $f:(Y,N)\rightarrow (X,M)$ is a $G$-torsor for the Kummer log flat topology, and since $T\mapsto \Gamma(T,\mathcal{O}_T)$ is a sheaf for the Kummer log flat topology \cite[\S{}3]{kato2}, the theory of descent along torsors gives us the relation $(f_*\mathcal{O}_Y)^G=\mathcal{O}_X$, hence the result.

$(ii)$ Let $x\in X$, and let $y\in Y$ be a point above $x$. Assume $x$ and $y$ have the same residue field, that we denote by $k$. Then the inertia group $I(y)$ of the action at $y$ is finite flat over $k$, and, according to Proposition~\ref{vistoli_prop}, becomes diagonalizable over a finite extension of $k$. Since being diagonalizable locally for the fppf topology and the \'etale topology are the same, $I(y)$ becomes diagonalizable over a finite separable extension of $k$. It follows that there exists an \'etale morphism $W\rightarrow X$ having $x$ in its image and a diagonalizable group scheme $\Delta$ over $W$ that extends $I(y)$. Therefore, according to the proof of \cite[Prop.~3.6]{aov}, there exists a finite scheme $T\rightarrow W$ on which $\Delta$ acts, together with an isomorphism $[Y/G]\times_X W\simeq [T/\Delta]$ of algebraic stacks over $W$. In the general case, the stack $[Y/G]$ being tame, it is possible, by \cite[Prop.~3.7]{aov}, to find a point $y\in Y$ above $x$ that is defined over a finite separable extension of the residue field of $x$. After performing an \'etale base change having $x$ in its image, we recover the situation where $x$ and $y$ have the same residue field, which concludes the proof.
\end{proof}

\begin{remark}
As we said before, Corollary~\ref{vistoli_cor} $(ii)$ is a generalization of Abhyankar's Lemma. In the version given by Grothendieck and Murre (see \cite[Thm.~2.3.2]{gm71}), the action of $G$ on $Y$ is locally induced (for the \'etale topology) from some action of a diagonalizable group $\Delta$. In our version, the stack is locally isomorphic (for the \'etale topology) to some $[T/\Delta]$, which is a slightly weaker statement. Let us also note that Chinburg, Erez, Pappas and Taylor prove in \cite[Thm.~6.4]{cept} another generalization of Abhyankar's Lemma, that they call a slice theorem, but they work with the fppf topology and assume that $G$ is commutative.
\end{remark}


\section{Tame covers with respect to a divisor}
\label{tame_divisor_section}


\subsection{Proof of Theorem~\ref{thm:main_intro}}

In this subsection, we consider again the situation of subsection \ref{XlogD_subsection} where $X$ is a noetherian regular scheme, and $D$ is a normal crossing divisor on $X$. As usual, $G$ is a finite flat group scheme over $X$.

\begin{definition}
\label{tame_cover_def}
A tame $G$-cover relative to $D$ is a finite flat scheme $Y\rightarrow X$ together with an action of $G$ on $Y$, such that
\begin{enumerate}
\item[$(1)$] the quotient stack $[Y/G]$ is tame, with moduli space $X$;
\item[$(2)$] $Y_U\rightarrow U$ is a $G_U$-torsor for the fppf topology.
\end{enumerate}
\end{definition}

\begin{remark}
\label{sophie_remark}
\begin{enumerate}
\item It has been proved by Marques \cite{sophie} that, given a finite flat morphism $Y\rightarrow X$ together with an action of $G$ on $Y$ such that $X=Y/G$, the quotient stack $[Y/G]$ is tame (with moduli space $X$) if and only if the action of $G$ on $Y$ is CEPT-tame \cite[Def.~7.1]{cept}. Therefore, we could replace condition $(1)$ in Definition~\ref{tame_cover_def} above by: $X=Y/G$ and the action of $G$ on $Y$ is CEPT-tame.
\item When $G$ is a finite commutative group scheme defined over a field, a notion of \emph{tamely ramified $G$-torsor} has been introduced by Biswas and Borne  \cite[Definition~2.2]{BiswasBorne2022}. By definition, a tamely ramified torsor is fppf locally (on the base) induced by some Kummer cover, hence the corresponding quotient stack is tame. So, tamely ramified $G$-torsors are tame $G$-covers in the sense of Definition~\ref{tame_cover_def}.
\end{enumerate}
\end{remark}

Of course, tame $G$-covers relative to $D$ together with $G$-equivariant morphisms form a category, that we denote by $\TC(X,D,G)$.

We have two restriction functors
$$
\begin{CD}
R_U:\Tors_{\kpl}(X(\log D),G) @>>> \Tors_{\fppf}(U,G) \\
\end{CD}
$$
and
$$
\begin{CD}
R'_U:\TC(X,D,G) @>>> \Tors_{\fppf}(U,G). \\
\end{CD}
$$

It follows from Theorem~\ref{thm_torsrep}, Proposition~\ref{prop:flatness} and Corollary~\ref{vistoli_cor} that the underlying scheme of a Kummer log flat $G$-torsor over $X(\log D)$ is a tame $G$-cover relative to $D$. Thus, we obtain a functor
$$
\begin{CD}
\Omega:\Tors_{\kpl}(X(\log D),G) @>>> \TC(X,D,G). \\
\end{CD}
$$
We immediately observe that
$$
R_U= R'_U\circ \Omega.
$$

\begin{proposition}
\label{prop:fullfaith}
The functors $R_U$ and $\Omega$ are fully faithful.
\end{proposition}

\begin{proof}
In order to prove that $R_U$ is fully faithful, it suffices to prove that the map $G(X)\to G(U)$ is bijective. This map is injective because $G\to X$ is affine, hence separated, and $U$ is dense in $X$. The surjectivity can be proved as follows: the base scheme $X$ is regular, $U$ is dense in $X$, and the morphism $G\to X$ is finite, so any section $U\to G$ can be extended into a section $X\to G$ according to \cite[6.1.14]{EGA2}.

Finally, it follows from the relation $R_U= R'_U\circ \Omega$ that $\Omega$ is faithful. Let us now prove the fullness of $\Omega$. Let $T$ and $T'$ be  be two objects in the category $\Tors_{\kpl}(X(\log D),G)$, and let $\psi:\Omega(T)\to \Omega(T')$ be a morphism in the category $\TC(X,D,G)$. We need to show that $\psi$ comes from a morphism in $\Tors_{\kpl}(X(\log D),G)$. We observe that $\psi_U:\Omega(T)_U\to \Omega(T')_U$ is a morphism between $G_U$-torsors, hence is an isomorphism. The functor $R_U$ being fully faithful, there exists a unique morphism $\phi:T\to T'$ in the category $\Tors_{\kpl}(X(\log D),G)$ extending $\psi_U$. Thus, we obtain two morphisms $\psi$ and $\Omega(\phi)$ in the category $\TC(X,D,G)$, whose restriction to $U$ is the same. We claim that $\psi=\Omega(\phi)$. The question is local on $X$, so we may assume that $X=\Spec(A)$ is the spectrum of an integral regular ring, that $U=\Spec(A[1/h])$ is a distinguished open set, and that $\Omega(T)=\Spec(B)$ and $\Omega(T')=\Spec(B')$ are spectra of finite rank free $A$-algebras. Then since $A\hookrightarrow A[1/h]$ the same holds for $B$ and $B'$, more precisely $B\hookrightarrow B\otimes_A A[1/h]$ and $B'\hookrightarrow B'\otimes_A A[1/h]$. By construction, the $A$-algebra morphisms $B'\to B$ corresponding to $\psi$ and $\Omega(\phi)$ agree after tensoring by $A[1/h]$, hence they are equal. In other terms $\psi=\Omega(\phi)$, which finishes the proof.
\end{proof}

The question now is: what is the image of this functor $\Omega$, that is, what kind of tame covers are coming from Kummer log flat torsors? In fact, the situation is very similar to what happens in the classical theory of tame covers, namely, Kummer log flat torsors provide initial objects in the fibers of the functor $R'_U$.

\begin{theorem}
\label{main_thm}
The restriction functors $R_U$ and $R'_U$ in the diagram
 \begin{equation*}
  \xymatrix@1@C+3mm{\Tors_{\kpl}(X(\log D),G) \ar[d]_(0.5){\Omega}\ar[rd]^(0.5){R_U} \\
  \TC(X,D,G) \ar[r]_(0.56){R'_U} & \Tors_{\fppf}(U,G)}
 \end{equation*}
have the same image. Moreover, given $t\in\Tors_{\fppf}(U,G)$ in this image, and letting $t^{\log}$ be the unique Kummer log flat $G$-torsor extending $t$, the  tame $G$-cover $\Omega(t^{\log})$ underlying $t^{\log}$ is the initial object in the category of tame $G$-covers extending $t$.
\end{theorem}

\begin{proof}
Everything follows from Lemma~\ref{lemma2}.
\end{proof}

\begin{remark}
Let us underline the fact that $R'_U$ is not fully faithful: given a torsor $t$ that lifts into a tame cover, there exists a priori several tame covers lifting $t$. According to Theorem~\ref{main_thm}, the Kummer log flat torsor furnishes an initial object in the category of all possible lifts. Thus, the theory of Kummer log flat torsors gives a canonical approach to the problem of extending torsors over $U$ into tame covers of $X$.
\end{remark}

We first prove Theorem~\ref{main_thm} in the case when $G$ is linearly reductive.

\begin{lemma}
\label{lemme_clef}
Let $G$ be a finite flat linearly reductive group scheme on $X$, and let $Y\rightarrow X$ be a tame $G$-cover relative to $D$. Let $Y_U^{\sharp}\rightarrow X(\log D)$ be the unique $G$-torsor for the Kummer log flat topology extending $Y_U\rightarrow U$ (whose existence is granted by Theorem~\ref{equivalence1}), and let $\widetilde{Y_U}$ be the underlying scheme of $Y_U^{\sharp}$. Then there exists a unique $G$-equivariant morphism $\widetilde{Y_U}\rightarrow Y$ whose restriction to $U$ is the identity.
\end{lemma}

\begin{proof}
The schemes $\widetilde{Y_U}$ and $Y$ are both finite flat (and in particular affine) over $X$, and the category of affine schemes over $X$ is a stack for the \'etale topology, hence the question is local for the \'etale topology on $X$. So, we may (and do) assume that $X=\Spec(R)$ is the spectrum of a strictly henselian local ring $R$. In this case, $Y=\Spec(A)$ and $\widetilde{Y_U}=\Spec(B)$ where $A$ and $B$ are finite flat $R$-algebras, hence are free $R$-algebras of finite rank. By flatness, $A$ and $B$ can be considered as $R$-subalgebras of $A_U=B_U$. We claim that  it suffices to prove   that $A\subseteq B$ in order to prove the statement. Indeed, if there exists a morphism $A\to B$ extending the identity $A_U=B_U$, then this morphism is unique (by flatness), hence corresponds to the inclusion $A\subseteq B$. We claim that the corresponding morphism of schemes $\psi:\Spec(B)\to\Spec(A)$ is $G$-equivariant: consider the two maps $G\times_X \widetilde{Y_U} \to Y$ given by $(g,y)\mapsto g\psi(y)$ and $(g,y)\mapsto \psi(gy)$. These maps coincide on $U$, and the schemes involved are finite flat over $X$, hence the maps coincide everywhere by the same argument as in the proof of Proposition~\ref{prop:fullfaith}.

Note that $R$ is a noetherian regular ring, because $X$ is noetherian regular.  Hence $R=\bigcap_{\mathfrak{p}\in \Sigma} R_\mathfrak{p}$, where $\Sigma$ denotes the subset of $X$ consisting of prime ideals of codimension $1$. If $\mathfrak{p}$ is such a prime, we let $A_\mathfrak{p}:=A\otimes_R R_\mathfrak{p}$. Since $A$ is a free $R$-module of finite rank, we have
$$
A = \bigcap_{\mathfrak{p}\in\Sigma}A_\mathfrak{p},
$$
and similarly for $B$. Thus we are reduced to prove that $A_\mathfrak{p}\subseteq B_\mathfrak{p}$ for all points $\mathfrak{p}$ of codimension $1$ in $D$. That is, we may assume that $X=\Spec(R)$ is the spectrum of a strictly henselian discrete valuation ring endowed with the canonical log structure.

On the other hand, the group $G$ being linearly reductive, we have an exact sequence
$$
\begin{CD}
1 @>>> \Delta @>>> G @>>> H @>>> 1 \\
\end{CD}
$$
where $\Delta$ is diagonalizable, and $H$ is constant of order coprime to the residue characteristic of $X$. Let $\mathfrak{X}=\prod_{i=1}^s\Z/n_i\Z$ be the character group of $\Delta$, i.e. $\Delta=\prod_{i=1}^s\mu_{n_i}$. According to \cite[expos\'e VIII, \S{}4]{SGA3-2}, the action of $\Delta$ on $Y$ corresponds to an $\mathfrak{X}$-grading of $A$, namely
\[
A=\bigoplus_{\underline{a}=(a_1,\dots,a_s)\in\mathfrak{X}}E_{\underline{a}},
\]
where $E_0$ is a subring of $A$ and the $E_{\underline{a}}$ are $E_0$-modules. As a direct summand of $A$, each $E_{\underline{a}}$ is finite locally free over $R$, hence is free of finite rank. In particular, $E_0$ is free of finite rank over $R$. Note that the quotient of $Y$ by the action of $\Delta$ (in the sense of Definition~\ref{def:tamestack}) is the affine scheme $Y/\Delta=\Spec(E_0)$.

On the other hand, we denote by $Y_U^{\sharp}/\Delta$ the quotient of $Y_U^{\sharp}$ by the action of $\Delta$ in the category of sheaves for the Kummer log flat topology. By Galois theory, $Y_U^{\sharp}/\Delta\to X(\log D)$ is a Kummer log flat $H$-torsor, hence (Theorem~\ref{thm_torsrep}) is representable by a log scheme.

We now consider the following commutative diagram in the category of log schemes (where the right-hand side consists of schemes with trivial log structure), and we want to prove the existence of a map on the top. 
$$
\begin{CD}
Y_U^{\sharp} @. Y \\
@VVV @VVgV \\
Y_U^{\sharp}/\Delta @>\nu_\Delta>> Y/\Delta \\
@VVV @VVhV \\
X(\log D) @>>> X \\
\end{CD}
$$

The existence and unicity of a map $\nu_\Delta$ in the middle of this diagram whose restriction to $U$ is the identity can be proved as follows: $Y_U^{\sharp}/\Delta\to X(\log D)$ is a Kummer log flat $H$-torsor, and $H$ is constant of order coprime to the residue characteristic of $X$, hence this torsor is Kummer log \'etale and in particular (by Lemma~\ref{log_etale_torsor}) the underlying scheme of $Y_U^{\sharp}/\Delta$ is the normalization of $X$ in $Y_U/\Delta$. Since $Y/\Delta\to X$ is finite flat, hence integral, the map $\nu_\Delta$ exists, and is unique, by the universal property of the normalization.

Let $\Spec(C)$ be the underlying scheme of $Y_U^{\sharp}/\Delta$, then $C$ is a normal $R$-algebra of finite rank, hence it can be written as product of discrete valuation rings. For simplicity we may, and do, assume that $C$ is a discrete valuation ring, in which case $Y_U^{\sharp}/\Delta$ is the spectrum of $C$ endowed with the canonical log structure. More precisely, $C$ is the ring of integers of a finite (tame) extension $L$ of the fraction field $K$ of $R$. By construction $L$ is the fraction field of $E_0$ hence we have $E_0\subseteq C$ and $C\otimes_R K = E_0\otimes_R K=L$. The $R$-module $A$ being of finite rank, it can be embedded in $A\otimes_R K$, and the action of $E_0$ on $A$ induces a structure of $L$-vector space on $A\otimes_R K$. Let $C\cdot A$ be the smallest
subset of $A\otimes_R K$ containing $A$ and stable by multiplication by $C$, where $C$ is seen as a subring of $L$. Then $C\cdot A$ is an $\mathfrak{X}$-graded $C$-algebra containing $A$, more precisely it is the image of the natural map $A\otimes_{E_0} C \to A\otimes_R K$.

Since a vector space has no torsion, $C\cdot A$ has no $C$-torsion. Since $C$ is a discrete valuation ring, $C\cdot A$ is flat as a $C$-module, hence is free of finite rank.

In order to prove that $A\subseteq B$ (as $R$-algebras), it now suffices to show that $C\cdot A \subseteq B$ as $C$-algebras.
Since $C\cdot A$ is an $\mathfrak{X}$-graded algebra over $C$, we are reduced to the case when $C=R$ and $G=\Delta$.

Let $\varepsilon:({\rm fs}/X)_{\kpl}\to ({\rm fs}/X)_{\fppf}$ be the forgetful map of sites. Let us choose a uniformizing element $\pi$ of $R$, then $\pi^\N\hookrightarrow R$ is a chart of the log structure of $X$. By \cite[App. D, Prop. D.1]{w-z1}, the chart gives a decomposition
\[
H^1_{\kpl}(X,\Delta)\cong H^1_{\fppf}(X,\Delta)\oplus (\oplus_{i=1}^s\pi^{\Z/n_i\Z}),
\]
under which the torsor $Y_U^{\sharp}\to X(\log D)$ is the contracted product of a classical fppf $\Delta$-torsor $T^{\mathrm{cl}}$ and a Kummer log flat $\Delta$-torsor $X(\log D)_{\underline{r}}$ corresponding to some $\underline{r}=(r_i)\in \oplus_{i=1}^s{\Z/n_i\Z}$. Let 
\[R[\underline{t}]/(\underline{t}^{\underline{n}}-\pi^{\underline{r}}):=R[t_1,\dots,t_s]/(t_1^{n_1}-\pi^{r_1},\dots,t_1^{n_s}-\pi^{r_s}).\]
Then the torsor $X(\log D)_{\underline{r}}$ is given by the saturation of the fine log scheme $\Spec (R[\underline{t}]/(\underline{t}^{\underline{n}}-\pi^{\underline{r}}))$ endowed with the log structure associated to $\langle t_1,\dots,t_s,\pi\rangle\hookrightarrow R[\underline{t}]/(\underline{t}^{\underline{n}}-\pi^{\underline{r}})$.
Since  $Y=T^{\mathrm{cl}}\wedge^{\Delta}((T^{\mathrm{cl}})^{-1}\wedge^{\Delta}Y)$, we are further reduced to the case where $Y_U^{\sharp}=X(\log D)_{\underline{r}}$ by replacing $Y$ by $(T^{\mathrm{cl}})^{-1}\wedge^{\Delta}Y$.

In this case, $Y_U^{\sharp}=X(\log D)_{\underline{r}}$ and we compute the affine ring $B$ of $Y_U^{\sharp}$ as follows. Let $P$ be the monoid generated by the monoid $\pi^\N$ and symbols $t_1,\dots,t_s$ subject to the relations $t_i^{n_i}=\pi^{r_i}$, and let $n=\prod_{i=1}^sn_i$. We claim that the saturation $P^{\rm sat}$ of $P$ is equal to $Q:=\{x\in P^\mathrm{gp} \mid x^{n}\in \pi^\N\}$. The inclusion $Q\subset P^{\rm sat}$ is clear. Conversely, any element of $P^{\rm sat}$ can be written as $a/b$ with $a,b\in P$ and there exists some $m$ such that $(a/b)^m\in P$. By construction, the $n$-power of any element of $P$ lies in $\pi^{\N}$. Hence $a^n,b^n\in\pi^{\N}$, and thus $(a^n/b^n)^m=((a/b)^m)^n\in \pi^{\N}$. Since $\pi^{\N}$ is saturated, we must have $(a/b)^n\in \pi^{\N}$. This finishes the proof of the claim.

We have that $B=R[\underline{t}]/(\underline{t}^{\underline{n}}-\pi^{\underline{r}})\otimes_{\Z[P]}\Z[Q]$. Since the $\mathfrak{X}$-grading of $A$ is induced by that of $A_K=B_K=K[\underline{t}]/(\underline{t}^{\underline{n}}-\pi^{\underline{r}})$, we have
\[A=R\oplus (\bigoplus_{\underline{a}\in\mathfrak{X}\backslash\{0\}}A\cap K\cdot \underline{t}^{\underline{a}}).\]
Thus it suffices to show that $A\cap K\cdot \underline{t}^{\underline{a}}\subset B$ for each $\underline{a}\in\mathfrak{X}$. Any element of $A\cap K\cdot \underline{t}^{\underline{a}}$ can be written as $u\underline{t}^{\underline{a}}\pi^j$ with $u\in R^\times$ and $j\in\Z$. We may assume that $u=1$. Since we must have $(\underline{t}^{\underline{a}}\pi^j)^n\in R$, we get $(\underline{t}^{\underline{a}}\pi^j)^n\in\pi^\N$, \emph{i.e.} $\underline{t}^{\underline{a}}\pi^j\in Q$. This finishes the proof.
\end{proof}

We now drop the assumption that $G$ is linearly reductive. The following Lemma is a reformulation of Theorem~\ref{main_thm}.

\begin{lemma}
\label{lemma2}
Let $Y\rightarrow X$ be a tame $G$-cover relative to $D$. Then there exists a $G$-torsor $Y_U^{\sharp}\rightarrow X(\log D)$ for the Kummer log flat topology that extends $Y_U\rightarrow U$. Moreover, if we denote by $\widetilde{Y_U}\rightarrow X$ the tame cover underlying this torsor, then there exists a unique $G$-equivariant map of tame covers $\widetilde{Y_U}\rightarrow Y$ whose restriction to $U$ is the identity.
\end{lemma}

\begin{proof}
Since $R_U$ is fully faithful (Proposition~\ref{prop:fullfaith}), the torsor $Y_U^{\sharp}$, if it exists, is unique up to unique isomorphism, hence the question is local for the \'etale topology on $X$. Therefore, the stack $[Y/G]$ being tame, we may assume, by \cite[Theorem~3.2]{aov}, that there exists a scheme $Z$, finite over $X$, and a linearly reductive group scheme $H$ acting on $Z$ such that $[Y/G]\simeq [Z/H]$. Let us denote by $M$ this stack, then we have a canonical $G$-torsor $Y\rightarrow M$ and a canonical $H$-torsor $Z\rightarrow M$. If we restrict everything to $U\subseteq X$, then, $Y_U\rightarrow U$ being a $G$-torsor, we get that $[Y_U/G]=U$, in other words $M_U=U$; in particular, $Z_U\rightarrow U$ is an $H$-torsor, obtained by restriction of the $H$-torsor $Z\to M$. We denote by $Z_U^{\sharp}\rightarrow X(\log D)$ the unique Kummer log flat $H$-torsor extending $Z_U\rightarrow U$, whose existence is granted by Theorem~\ref{equivalence1}.

Consider the fiber product $E:=Y\times_M Z$, then by construction the projection $E\to Y$ is an $H$-torsor and the projection $E\to Z$ is $G$-torsor, in particular $E$ is a scheme. Moreover, we observe that $E\to M$ is a $G\times H$-torsor, since it is the fiber product of a $G$-torsor and an $H$-torsor over the base stack $M$.

By construction the map $Z\rightarrow X$ is finite, we claim that it is also flat. Indeed, we have a commutative diagram
$$
\begin{CD}
E @>>> Z \\
@VVV @VVV \\
Y @>>> X \\
\end{CD}
$$
in which $E\rightarrow X$ is finite flat (because $E\rightarrow Y$ and $Y\rightarrow X$ are finite flat), and $E\rightarrow Z$ is also finite flat, hence is faithfully flat. By \cite[2.2.13~(iii)]{EGA4.2} we deduce that $Z\rightarrow X$ is flat. Therefore, $Z\rightarrow X$ is a tame $H$-cover relative to $D$, and, according to Lemma~\ref{lemme_clef}, there exists a unique morphism $Z_U^{\sharp}\rightarrow Z$ whose restriction to $U$ is the identity.

Let $E_U^{\sharp}$ be the fiber product $E\times_Z Z_U^{\sharp}$, computed in the category of fs log schemes. Then the projection $E_U^{\sharp}\rightarrow Z_U^{\sharp}$ is a strict morphism, and is a Kummer log flat $G$-torsor, obtained by base-change of the fppf $G$-torsor $E\rightarrow Z$. Furthermore, since $H$ acts equivariantly on $E\to Z$ and $Z_U^{\sharp}\rightarrow Z$, it acts diagonally on their fiber product $E_U^{\sharp}$.

The situation is now the following: $G\times H$ acts on $E_U^{\sharp}$, the map $E_U^{\sharp}\rightarrow Z_U^{\sharp}$ is a $G$-torsor, and the map $Z_U^{\sharp}\rightarrow X(\log D)$ is an $H$-torsor. We now want to show that $E_U^{\sharp}\rightarrow X(\log D)$ is a $G\times H$-torsor. Let us pick a Kummer log flat cover  $T\rightarrow X(\log D)$ together with a section $s:T\rightarrow E_U^{\sharp}$ (we could take $T=E_U^{\sharp}$ for example). Then this section $s$ induces a section $s':T\rightarrow Z_U^{\sharp}$. Therefore, $Z_U^{\sharp}\times T$ is the trivial $H$-torsor over $T$, \emph{i.e.} $s'$ induces an isomorphism $Z_U^{\sharp}\times T\simeq H\times T$ (the unadorned fiber products  are over $X(\log D)$). Also, the section $s$ together with the action of $H$ on $E_U^{\sharp}$ induces a section $s_H:H\times T\rightarrow E_U^{\sharp}\times T$. By composing this $s_H$ with the previous isomorphism, we conclude that $E_U^{\sharp}\times T\rightarrow Z_U^{\sharp}\times T$ is the trivial $G$-torsor. Putting all that together, we find that the natural map
$$
G\times H\times T\rightarrow  E_U^{\sharp}\times T
$$
induced by the section $s$ is an isomorphism. Thus, by Kummer log flat descent, this proves that $E_U^{\sharp}\rightarrow X(\log D)$ is a $G\times H$-torsor.

To conclude, as in any topos, a $G\times H$-torsor can be written uniquely as the product of a $G$-torsor and an $H$-torsor. Thus, there exists a $G$-torsor $Y_U^{\sharp}\to X(\log D)$ such that
$$
E_U^{\sharp} = Y_U^{\sharp} \times_{X(\log D)} Z_U^{\sharp}.
$$

Finally, the natural $H$-equivariant map $E_U^{\sharp}\rightarrow E$ gives rise, after modding out by $H$ on both sides, to a map $Y_U^{\sharp}\rightarrow Y$ whose restriction to $U$ is the identity, which concludes the proof.
\end{proof}

In the statement of Lemma~\ref{lemma2}, the notation suggests that $\widetilde{Y_U}$ is the normalization of $X$ in $Y_U$. This is indeed the case when $G$ is constant (Lemma~\ref{log_etale_torsor}). At the other end of the spectrum, if $G=\mu_n$ for some $n$ which is not invertible on $X$, then this no longer holds in general, as show the examples below.

\begin{example}
Let $p$ be an odd prime number and let $K$ be a finite extension of $\Q_p$ which contains a primitive $p^2$-th root of unity. Let $L/K$ be a totally ramified Kummer extension of degree $p^2$ (which corresponds to a $\mu_{p^2}$-torsor by classical Kummer theory), and let $\mathcal{O}_L$ and $\mathcal{O}_K$ be the respective rings of integers. In \cite{miyata97}, the Kummer order $\widetilde{\mathcal{O}_L}$ is defined as being the largest subring of $\mathcal{O}_L$ which is $(\Z/p^2\Z)$-graded. Then $\Spec(\widetilde{\mathcal{O}_L})\to \Spec(\mathcal{O}_K)$ is the underlying scheme morphism of the Kummer log flat $\mu_{p^2}$-torsor extending $\Spec(L)\to \Spec(K)$, and, according to \cite[Theorem~2]{miyata97}, the equality $\widetilde{\mathcal{O}_L} = \mathcal{O}_L$ does not hold in general.
\end{example}

\begin{example}
Let $R$ be a discrete valuation ring with fraction field $K$ of characteristic $p>0$, and a chosen uniformizer $\pi$. We endow $A=R[t]/(t^{p^2}-\pi^p)$ with the canonical $\mu_{p^2}$-action. We note that $A$ is a finite flat $R$-algebra, that the generic fiber of $\Spec(A)\to \Spec(R)$ is a $\mu_{p^2}$-torsor, and that $\mu_{p^2}$ is linearly reductive; so $\Spec(A)$ is a tame $\mu_{p^2}$-cover of $\Spec(R)$ relative to the closed point of $\Spec(R)$. Apparently $t^p-\pi$ is nilpotent in  $A$, so for any $a\in K$ the element $a(t^p-\pi)$ is integral over $R$, in particular $\frac{t^p-\pi}{\pi^2}$ is integral over $R$. Let $B$ be the underling ring of the Kummer log flat $\mu_{p^2}$-torsor associated to $\Spec(A)$, then we have $B=A[\frac{t^p}{\pi}]\cong R[t,v]/(v^p-1,t^p-v\pi)$. Since $\frac{t^p-\pi}{\pi^2}\notin B$, we have that $B$ is strictly contained in the integral closure of $A$ inside $A\otimes_RK=K[t]/(t^{p^2}-\pi^p)$.
\end{example}


\subsection{Tame objects for Hopf algebras over Dedekind rings}
\label{CH_tameries}

Let $R$ be a Dedekind ring with fraction field $K$, and let $X=\Spec(R)$. Let $G$ be a finite flat group scheme over $X$, and let $H$ be the $R$-Hopf algebra of $G$, so that $G=\Spec(H)$. We assume that the generic fiber of $G$ is an \'etale group scheme over $\Spec(K)$ (this is satisfied if $K$ has characteristic $0$).

The notion of Galois $H$-object and tame $H$-object we refer to in this section are those defined in the paper by Childs and Hurley \cite{ch86}.

Let $\Spec(F)\rightarrow \Spec(K)$ be a $G$-torsor over $\Spec(K)$. In other terms, $F$ is a Galois $H_K$-object. Many authors have been interested in the following question: under which condition is it possible to find an $R$-order in the algebra $F$ that has the structure of tame $H$-object?

Let $Z\rightarrow X$ be the normalization of $X$ in $\Spec(F)$ (\emph{i.e.} the spectrum of the integral closure of $R$ in $F$). Let $U$ be the largest open subset of $X$ such that $G_U\rightarrow U$ and $Z_U\rightarrow U$ are \'etale. We note that $U$ is dense in $X$, because these conditions are satisfied at the generic point of $X$. Therefore, the set $D:=X-U$ is a finite set of prime ideals of $R$, and is trivially a divisor with normal crossings on $X$. Moreover, $Z_U\rightarrow U$ has a natural structure of $G_U$-torsor: this is clear when $G_U$ is constant. In the general case we perform an \'etale base change which turns $G_U$ into a constant group, and conclude by descent theory.

We can now state the following consequence of Theorem~\ref{main_thm}:

\begin{proposition}
\label{CH_prop}
It is possible to find an $R$-order in the algebra $F$ that has the structure of tame $H$-object if and only if it is possible to extend $Z_U\rightarrow U$ into a $G$-torsor over $X(\log D)$ for the Kummer log flat topology. Moreover, the algebra underlying the Kummer log flat torsor is the maximal $R$-order in $F$ that has a structure of tame $H$-object.
\end{proposition}

\begin{proof}
By combining \cite[Thm.~2.6]{cept} and the result by S. Marques (see \cite{sophie} and Remark~\ref{sophie_remark}), we see that an $R$-order $\mathcal{O}$ in $F$ is a tame $H$-object if and only if $\Spec(\mathcal{O})\rightarrow X$ is a tame $G$-cover relative to $D$. The result then follows from Theorem~\ref{main_thm}.
\end{proof}

\begin{remark}
In a different direction, one may ask the following: when is it possible to find some Hopf order $H'$ in $H_K$ such that the integral closure of $R$ in $F$ is a tame $H'$-object? In our language: when is it possible to find a finite flat group scheme $G'$ over $X$ extending $G_U$ such that $Z\rightarrow X$ is a tame $G'$-cover?
This question has a complete known answer in special cases, for example, when the generic fiber of $G$ is constant cyclic of prime order (see \cite{childs87} and \cite{byott95}).
The question we have been investigating here goes the other way: we fix $H$, and ask when $F$ can be extended into a tame $H$-object.
\end{remark}


\subsection{Tame covers of arithmetic schemes}
\label{CEPT_tameries}

Let $R$ be a Dedekind ring with fraction field $K$, and let $X\rightarrow \Spec(R)$ be a regular scheme, projective and flat of relative dimension $d$ over $\Spec(R)$. We assume that the generic fiber of $X$ is smooth over $\Spec(K)$ (this is satisfied if $K$ has characteristic $0$).

Let $G$ be a finite flat group scheme over $\Spec(R)$, and let $\pi:Y\rightarrow X$ be a finite morphism with an action of $G$ on $Y$ that is CEPT-tame, with quotient $X$.

We make the following hypotheses:

\begin{enumerate}
\item[(Hyp~1).] The morphism $\pi$ is flat.
\item[(Hyp~2).] There exists a normal crossing divisor $D$ on $X$, with complement $U$, such that $\pi_U:Y_U\rightarrow U$ is a $G$-torsor.
\end{enumerate}

We note that similar hypotheses have been considered in several papers by Cassou-Nogu\`es, Chinburg, Pappas and Taylor \cite{cpt07}, \cite{cpt09}, \cite{ct10}. In most cases, they also assume that $D$ is a fibral divisor. Let us note that when $\Spec(R)\rightarrow \Spec(\Z)$ is surjective, then $Y_K\rightarrow X_K$ is a $G$-torsor (see \cite[1.2.4(d)]{cept97}), so $D$ is fibral.

If $G$ is constant and $Y$ is normal, then $\pi$ is flat (see \cite[Remark~3.1]{cpt09}). In this case, $Y\rightarrow X$ is tame in the sense of Grothendieck and Murre, so there exists a log structure on $Y$ that turns it into a Kummer log flat torsor over $X(\log D)$. The following Proposition is a generalization of this result. Namely, starting from a CEPT-tame $G$-cover satisfying (Hyp~1) and (Hyp~2), there always exists a procedure (which may be thought of as a kind of normalization) which allows to build a Kummer log flat $G$-torsor from it. Moreover, the object obtained is initial in the category of CEPT-tame $G$-covers.

\begin{proposition}
\label{CEPT_prop}
Let $\pi:Y\rightarrow X$ be a CEPT-tame $G$-cover satisfying (Hyp~1) and (Hyp~2) above. Then there exists a unique $G$-torsor for the Kummer log flat topology over $X(\log D)$ that extends $\pi_U:Y_U\rightarrow U$. Moreover, the cover underlying this torsor is CEPT-tame, and is initial in the category of CEPT-tame $G$-covers extending $\pi_U$.
\end{proposition}

\begin{proof}
By the result of S. Marques (see \cite{sophie} and Remark~\ref{sophie_remark}), our tameness coincides with CEPT-tameness in the present setting. The result follows from Lemma~\ref{lemma2}.
\end{proof}

We hope that Proposition~\ref{CEPT_prop} will open new perspective on CEPT-tameness. In particular, the ``slice theorem'' (Corollary~\ref{vistoli_cor}) gives us a nice description of the local structure, for the \'etale topology, of the initial CEPT-tame object.

Assume now that $d=1$, that is, $X$ is an arithmetic surface. Then, starting from a CEPT-tame $G$-cover that is generically a torsor, and that satisfies (Hyp~1), it is possible to construct, up to blowing-up a finite number of closed points of $X$, a CEPT-tame $G$-cover satisfying (Hyp~1) and (Hyp~2).

More precisely, since $\pi:Y\rightarrow X$ is generically a torsor, there exists a dense open subset $V\subseteq X$ such that $\pi_V$ is a torsor. Since $X$ is a regular scheme of dimension $2$, it is not very hard to show that there exists a divisor $D$ on $X$ such that $X=V\cup D$. Let us note that, if the generic fiber of $G$ is \'etale, then we may take $D$ to be the branch divisor of $\pi$, namely the smallest divisor on $X$ outside which $\pi$ is \'etale.
If $D$ is not fibral, we also assume that $R$ is excellent. Then, according to \cite[Chap.~9, Thm.~2.26]{liu}, by blowing-up a finite number of suitably chosen closed points of $X$, one obtains a projective birational morphism $f:X'\rightarrow X$ where $X'$ an arithmetic surface, such that $f^*D$ is a divisor with normal crossings on $X'$, and $f$ is an isomorphism outside $f^*D$. Now, pulling-back $\pi$ along $f$, we get a CEPT-tame $G$-cover $\pi':Y'\rightarrow X'$ satisfying (Hyp~1) and (Hyp~2).



\bibliographystyle{alpha}
\bibliography{biblioLFT}



\bigskip

\textsc{Jean Gillibert}, Institut de Math{\'e}matiques de Toulouse, CNRS UMR 5219, 118 route de Narbonne, 31062 Toulouse Cedex 9, France.

\emph{E-mail address:} \texttt{jean.gillibert@math.univ-toulouse.fr}
\medskip

\textsc{Heer Zhao}, Institute for Advanced Study in Mathematics, Harbin Institute of Technology, Harbin, 150001, P. R. China.

\emph{E-mail address:} \texttt{logpt@sina.com}


\end{document}